\documentclass[12pt]{amsart}
\usepackage{amsfonts}
\usepackage{amssymb}
\usepackage{graphicx}

\def\la{\lambda}           
     \def\l{\lambda}

\def\D{{\mathbb D}}  
\def\C{{\mathbb C}}

\def\({\left(}       \def\){\right)}

\newtheorem{lem}{\sc Lemma}
\newtheorem{thm}{\sc Theorem}

\newtheorem{ex}{\sc Example}

\newenvironment{pf}{\noindent{\textit{Proof. }}}{$\Box$ }

\begin{document}
\title[Criteria for Univalence and quasiconformal extension]
{Criteria for univalence and quasiconformal extension of harmonic mappings in terms of the Schwarzian derivative}


\author[R. Hern\'andez]{Rodrigo Hern\'andez}
\address{Facultad de Ingenier\'{\i}a y Ciencias, Universidad Adolfo Ib\'a\~nez, Av. Padre Hurtado 750, Vi\~na del Mar, Chile.} \email{rodrigo.hernandez@uai.cl}

\author[M. J. Mart\'{\i}n]{Mar\'{\i}a J. Mart\'{\i}n}
\address{Department of Physics and Mathematics, University of Eastern Finland, P.O. Box 111, FI-80101 Joensuu, Finland.} \email{maria.martin@uef.fi}

\date{\today}
\thanks{This research was supported by grants Fondecyt $1110160$ and $1110321$, Chile. The second author is also supported by
Academy of Finland grant 268009 and by Spanish MINECO Research Project MTM2012-37436-C02-02.}

\begin{abstract}
We prove that if the Schwarzian norm of a given complex-valued locally univalent harmonic mapping $f$ in the unit disk is small enough, then $f$ is, indeed, globally univalent and can be extended to a quasiconformal mapping in the extended complex plane.
\end{abstract}
\maketitle


\section*{Introduction}
In 1949, Nehari \cite{Nehari} proved that if a locally univalent analytic function $\varphi$ in the unit disk $\D$ satisfies
\[
\sup_{z\in\D} |S(\varphi)(z)|\, (1-|z|^2)^2\leq 2\,,
\]
then $\varphi$ is globally univalent in $\D$. Here, $S(\varphi)$ denotes the \emph{Schwarzian derivative} of $\varphi$ defined by
\begin{equation}\label{eq-Schw}
S(\varphi)=\left(\frac{\varphi''}{\varphi'}\right)'-\frac 12 \left(\frac{\varphi''}{\varphi'}\right)^2\,.
\end{equation}
Ahlfors and Weill \cite{AW} generalized Nehari's criterion of univalence by proving that if such function $\varphi$ satisfies $\|S(\varphi)\|\leq 2t$ for some $t<1$, then $\varphi$ is injective in $\D$ and has a $K$-quasiconformal extension to $\widehat{\C}=\C\cup\{\infty\}$, where $K=(1+t)/(1-t)$.
\par\smallskip
Let now $f$ be a complex-valued locally univalent harmonic mapping in the unit disk. By considering the complex conjugate, if needed, we can assume that $f$ is sense-preserving. This is, $f=h+\overline g$ where $h$ and $g$ are analytic functions in $\D$ such that $h$ is locally univalent and the (second complex) dilatation $\omega=g'/h'$ is an analytic function mapping the unit disk into itself. The following definition for the Schwarzian derivative of such functions $f$ was presented in \cite{HM-Schw}:
\begin{equation}\label{eq-schharmonic}
S_f= S(h)+\frac{\overline \omega}{1-|\omega|^2} \left(\frac{h''}{h'}\, \omega'-\omega''\right)-\frac 32 \left(\frac{\overline \omega\, \omega'}{1-|\omega|^2}\right)^2\,,
\end{equation}
where $S(h)$ is the classical Schwarzian derivative of $h$, as in \eqref{eq-Schw}.
\par\smallskip
The main purpose in this note is to prove that there exists a constant $\delta_0 >0$ such that if the locally univalent harmonic mapping $f$ in the unit disk has Schwarzian norm
\[
\|S_f\|=\sup_{z\in\D} |S_f(z)|\, (1-|z|^2)^2\leq \delta_0\,,
\]
then $f$ is one-to-one in $\D$. We will also see that if $\|S_f\|\leq \delta_0 t$ for some $t<1$, then $f$ has a quasiconformal extension to $\widehat \C$.

\section{Background}\label{sec-background}

\subsection{The Schwarzian derivative}\label{ssec-Schwarzian} As it was mentioned in the introduction, every harmonic mapping $f$ in the unit disk $\D$ can be written as $f=h+\overline g$ with $h$ and $g$ analytic in $\D$. This decomposition is unique up to an additive constant (see \cite[p. 7]{Dur-Harm}). We refer to the reader to the book \cite{Dur-Harm} for a comprehensive treatment on harmonic mappings.
\par
Lewy \cite{Lewy} proved that a harmonic mapping in the unit disk is locally univalent if and only if its Jacobian is different from zero. In terms of the decomposition $f=h+\overline g$, the Jacobian $J_f$ of $f$ equals $|h'|^2-|g'|^2$. Thus, locally univalent harmonic mappings in $\D$ are either \emph{sense-preserving} if $J_f>0$ or \emph{sense-reversing} if $J_f<0$. Note that any analytic function is a sense-preserving harmonic mapping. Also, that a harmonic function $f=h+\overline g$ is sense preserving if and only if $h$ is locally univalent and the dilatation $\omega=g'/h'$ maps the unit disk into itself. It is obvious that $f$ is sense-preserving if and only if $\overline f$ is sense-reversing. In this paper, we will consider harmonic mappings which are sense-preserving in the unit disk. For this kind of mappings, the Schwarzian derivative is given by \eqref{eq-schharmonic}. It is clear that if $f$ is analytic, then $S_f$ coincides with the classical definition of the Schwarzian derivative given by \eqref{eq-Schw}.
\par
Several properties of this operator are the following.
\begin{enumerate}

\item[(i)] $S_f\equiv 0$ if and only if $f=\alpha T+\beta\overline{T}$, where $|\alpha|\neq|\beta|$ and $T$ is a M\"{o}bius transformation
$$T(z)=\frac{az+b}{cz+d}\,,\quad ad-bc\neq 0\,.$$

\item[(ii)] Whenever $f$ is a sense-preserving harmonic mapping and $\phi$ is an analytic function such that the composition $f\circ\phi$ is well-defined, the Schwarzian derivative of $f\circ\phi$ can be computed using the \emph{chain rule}
    $$S_{f\circ\phi}=S_f(\phi)\cdot(\phi')^2+S\phi\,.$$

\item[(iii)] For any \emph{affine mapping} $L(z)=az+b\overline z$ with $|a|\neq |b|$, we have that $S_{L\circ f}=S_f$. Note that $L$ is sense-preserving if and only if $|b|<|a|$.

\end{enumerate}

\par\smallskip
The \emph{Schwarzian norm} $\|S_f\|$ of a sense-preserving harmonic mapping $f$  in the unit disk is defined by
\[
\|S_f\|=\sup_{z\in\D} |S_f(z)|\cdot(1-|z|^2)^2\,.
\]
\par
It is easy to check (using the chain rule again and the Schwarz-Pick lemma) that $\|S_{f\circ\, \sigma}\|=\|S_f\|$ for any automorphism of the unit disk $\sigma$. For further properties of $S_f$ and the motivation for this definition, see \cite{HM-Schw}.
\subsection{An affine and linear invariant family}\label{ssec-ALIF}
In \cite{POM-I, POM-II} Pommerenke studied the so-called \emph{linear invariant
families\,}; that is, fa\-mi\-lies of locally univalent holomorphic functions $\varphi$ in the unit disk normalized by the conditions $\varphi(0)=1-\varphi^\prime(0)=0$ and
which are closed under the  transformation
$$\Phi_{\zeta}(z)=\frac{\varphi\left(\displaystyle\frac{\zeta+z}{1+\overline{\zeta}z}\right)
-\varphi(\zeta)}{(1-|\zeta|^2)\varphi^\prime(\zeta)}\,,\quad \zeta\in\D.$$
\par
Let $\mathcal{F}$ be a family of sense-preserving harmonic mappings $f=h+\overline{g}$ in $\D$, normalized with $h(0)=g(0)=0$ and $h'(0)=1$. The family is said to be {\it affine and linear invariant} if it closed under the two operations of \emph{Koebe transform} and \emph{affine change}:
\[
K_{\zeta}(f)(z)=\frac{f\left(\displaystyle\frac{z+\zeta}{1+\overline{\zeta} z}\right)-f(\zeta)}{(1-|\zeta|^2)h'(\zeta)}\,,\quad |\zeta|<1\,,
\]
and
\[
A_\varepsilon(f)(z)=\frac{f(z)-\overline{\varepsilon f(z)}}{1-\overline\varepsilon g'(0)}\,,\quad |\varepsilon|<1\,.
\]
We refer to the reader to the paper \cite{S-S} where Sheil-Small offers an in depth study of affine and linear invariant families $\mathcal{F}$ of harmonic mappings in $\D$.
\par
Using that the Schwarzian derivative for harmonic mappings satisfies the chain rule and is invariant under affine changes $af+b\overline{f}$, $|a|\neq |b|$, it is easy to show that the family $\mathcal{F}_{\lambda}$ of sense-preserving harmonic mappings $f=h+\overline{g}$ in $\D$ with $h(0)=g(0)=0, h'(0)=1,$ and $||S_f||\leq \lambda$ is affine and linear invariant. We let  $\mathcal F_\lambda ^0=\{f\in \mathcal F_\lambda\colon g^\prime(0)=0\}$. These two families are studied in \cite{CHM}, where it is shown in particular that $\mathcal F_\lambda ^0$ is a compact family of harmonic mappings with respect to the topology of uniform convergence on compact subsets of $\D$. In that paper, the following notation was used: an analytic function in the unit disk $\omega$ with $\omega(\D)\subset\D$ is said to belong to $\mathcal A_\lambda^0$ (\emph{resp.} $\mathcal A_\lambda$) if there exists a harmonic mapping $f=h+\overline g\in \mathcal F_\lambda ^0$ (\emph{resp.} $\mathcal F_\lambda$) with dilatation $\omega$. The quantity
\[
R_\lambda= \max_{\omega \in \mathcal A_\lambda^0}|\omega'(0)| = \sup_{\omega \in \mathcal A_\lambda} \|\omega^*\|\,,
\]
where
\[
\|\omega^*\|=\sup_{z\in\D} \frac{|\omega'(z)|\cdot(1-|z|^2)}{1-|\omega(z)|^2}\,,
\]
was shown to play a distinguished role in the analysis offered in \cite{CHM}.
\par
Perhaps at this point we should mention that according to the first of the properties for the Schwarzian derivative mentioned in the previous section, we have that the family $\mathcal F_0$ consists only of functions of the form $f=\alpha T+\beta\overline{T}$, where $|\alpha|>|\beta|$ and $T$ is a M\"{o}bius transformation. Therefore, all dilatations in $\mathcal A_0$ are constant functions.

\section{Main Results}\label{sec-mainresult}
The following lemma will be important for our purposes.
\begin{lem}\label{lem-Rlambda}
As before, let $R_\lambda=\max_{\omega \in \mathcal A_\lambda^0}|\omega'(0)|\,.$  Then
\begin{equation*}\label{eq-goal}
\lim_{\lambda \to 0^+} R_\lambda=0\,.
\end{equation*}
\end{lem}
\begin{pf}
Since $\mathcal F_{\lambda_1}^0\subset \mathcal F_{\lambda_2}^0$ whenever $0<\l_1 < \l_2$, we have $0\leq R_{\lambda_1} \leq R_{\lambda_2}$ as well. Therefore, we conclude that there exists $\lim_{\lambda \to 0^+} R_\lambda$ and it remains to check that this limit equals $0$.
\par
Consider an arbitrary positive number $\lambda$. By the definition of $R_\lambda$ and the compacity of $\mathcal F_\lambda^0$, we see that for each such $\lambda$ there is a harmonic mapping $f_\lambda\in \mathcal F_\lambda^0$ with dilatation $\omega_\lambda$ satisfying $|\omega'_\lambda(0)|=R_\lambda$. Since for a given $\rho>0$ the family $\{f_\lambda\colon  \lambda\leq \rho\} \subset \mathcal F_\rho^0$ and $\mathcal F_\rho^0$ is compact, we see that there is a function $f_0 \in \cap_{\rho >0} \mathcal F_\rho^0$ with dilatation $\omega_0$ such that $f_\lambda\to f_0$ as $\l\to 0$ uniformly on compact subsets in the unit disk (hence $\omega_\lambda'(0)\to \omega_0'(0)$ as $\l\to 0$ too). Obviously, $\cap_{\rho >0} \mathcal F_\rho ^0=\mathcal F_0^0$ and the dilatations of functions in $\mathcal \mathcal \mathcal F_0^0$ are constants, thus $0=\omega_0'(0)=\lim_{\lambda \to 0} R_\lambda$.
\end{pf}
\par\smallskip
We now state and prove the main theorems in this paper.
\begin{thm}\label{thm-main}
There exists $\delta_0 >0$ such that if $\|S_f\|\leq \delta_0$, then $f$ is univalent.
\end{thm}
\begin{pf}
For any real number $\lambda >0$, we have that if $f=h+\overline g\in \mathcal F_\lambda$, the Schwarzian norm of $h$ is bounded by \cite[Thm. 6]{HM-Schw}. Hence (see \cite{POM-I}),
\[
\sup_{z\in\D} \left|\frac{h''(z)}{h'(z)}\right|\, (1-|z|^2)\leq K_1
\]
for some constant $K_1>0$. Moreover, by using $\omega$ to denote the dilatation of $f$, we have (see, for instance, \cite{GZ}) that there exists another positive constant $K_2$ such that
\[
\sup_{z\in\D} \frac{\left|\omega''(z)\right|\, (1-|z|^2)^2}{1-|\omega(z)|^2}\leq K_2 \|\omega^*\|\leq K_2 R_\lambda\,.
\]
Hence, using \eqref{eq-schharmonic} and the triangle inequality, we see that for any such function $f=h+\overline g$,
\begin{equation}\label{eq-Sh estimate1}
\|S(h)\|\leq \lambda +K_1 R_\lambda +K_2 R_\lambda +\frac 32 R_\lambda^2\,.
\end{equation}
Now, using Lemma~\ref{lem-Rlambda} and the fact that $R_\la$ increases with $\lambda$, we have that there exists a unique solution $\delta_0$, say, of the equation
\[
\lambda +K_1 R_\lambda +K_2 R_\lambda +\frac 32 R_\lambda^2=2\,.
\]
This implies by \eqref{eq-Sh estimate1} that if $\lambda\leq \delta_0$, then $\|S(h)\|\leq 2$. In other words,
the analytic part $h$ of any function $f=h+\overline g\in \mathcal F_{\delta_0}$ is univalent by the classical Nehari criterion of univalence.
\par
To prove that not only $h$ but the function $f=h+\overline g$ itself is univalent whenever $f\in \mathcal F_{\delta_0}$, we proceed as follows. By the affine invariance property of $F_{\delta_0}$, we see that for any $a\in\D$ the function $f_a=f+\overline{af}$ belongs to $F_{\delta_0}$ as well. It is easy to check that if $f_a=h_a+\overline{g_a}$, then $h_a=h+\overline{a}g$. Thus, the functions $h+\overline{a}g$ are also univalent for all  $|a|<1$. Since $f$ is sense-preserving, an application of Hurwitz theorem gives that $h+\overline{a} g$ is indeed univalent for all $|a|\leq 1$. A direct application of \cite[Prop. 2.1]{HM-Schw} shows that the function $f=h+\overline g$ is univalent, as was to be shown.
\end{pf}
\par
We would like to point out that by finding an upper bound for the quantity $R_\lambda$ in terms of $\lambda$, one could give an estimate of the value $\delta_0$ in the previous theorem. Unfortunately, so far we are not able to obtain such upper bound.
\par
The next result is related to a criterion for quasiconformal extension of harmonic mappings in terms of their Schwarzian norm.

\begin{thm}\label{thm-main2}
Let $f$ be a sense-preserving harmonic mapping in the unit disk with $\|S_f\|\leq \delta_0 t$ for some $t<1$, where $\delta_0$ is as in Theorem~\ref{thm-main}. Assume, in addition, that the dilatation $\omega_f$ of $f$ satisfies
\begin{equation*}\label{eq-dilat}
\|\omega_f\|_\infty=\sup_{z\in\D}|\omega_f(z)|<1\,.
\end{equation*}
Then $f$ can be extended to a quasiconformal map in $\widehat\C$.
\end{thm}
Before proving this second theorem, we would like to stress that the hypotheses $\|\omega\|_\infty<1$ cannot be removed as the following example shows.
\begin{ex}
Consider the sense-preserving harmonic mapping $f=z+\overline g$, where $g^\prime$ equals the lens-map $\ell_\alpha$, $0<\alpha\leq 1$, defined by
\begin{equation*}\label{eq-lensmap}
\ell_\alpha(z)=\frac{\ell(z)^\alpha-1}{\ell(z)^\alpha+1}\,,\quad z\in\D\,,
\end{equation*}
with $\ell(z)=(1+z)/(1-z)$. Note that $\ell_1$ equals the identity in the unit disk and that $\|\ell_\alpha\|_\infty=1$ for all $0<\alpha\leq 1$. In \cite{HO}, it is explicitly checked that $\|\ell_\alpha^*\|=\alpha$.
\par\smallskip
Bearing in mind \eqref{eq-schharmonic}, that the dilatation of $f$ is $\ell_\alpha$, and that
\[
\sup_{z\in\D} \frac{\left|\ell_\alpha''(z)\right|\, (1-|z|^2)^2}{1-|\ell_\alpha(z)|^2}\leq K_2 \|\ell_\alpha^*\|= K_2 \alpha
\]
for some absolute constant $K_2$, we have
\begin{eqnarray*}
\|S_f\|&=&\sup_{z\in\D} \left|\frac{\overline{\ell_\alpha(z)}\, \ell_\alpha''(z)}{1-|\ell_\alpha(z)|^2}+\frac 32 \left(\frac{\overline{\ell_\alpha(z)}\, \ell_\alpha'(z)}{1-|\ell_\alpha(z)|^2}\right)^2\right|\, (1-|z|^2)^2\\
&\leq& \sup_{z\in\D} \frac{\left|\ell_\alpha''(z)\right|\, (1-|z|^2)^2}{1-|\ell_\alpha(z)|^2}+\frac 32 \sup_{z\in\D} \left|\frac{\ell_\alpha'(z)\, (1-|z|^2)}{1-|\ell_\alpha(z)|^2}\right|^2\\
&\leq& K_2 \alpha +\frac 32 \alpha^2\,.
\end{eqnarray*}
Therefore, by choosing any $\alpha$ small enough, we obtain $\|S_f\|\leq \delta_0 t$ for any given $0<t<1$. On the other hand, the function $f$ is not quasiconformal since its (sencond complex) dilatation coincides with $\ell_\alpha$ and $\|\ell_\alpha\|_\infty=1$.
\end{ex}
\par
We now prove Theorem~\ref{thm-main2}.
\par\smallskip

\begin{pf}
Since we are assuming that $\|S_f\|\leq \delta_0 t$ for some $t<1$, we have that $f\in \mathcal F_{\delta_0 t}$. By arguing as in the proof of the previous theorem and using again that if $\lambda_1\leq \lambda_2$ then $R_{\lambda_1}\leq R_{\lambda_2}$, we get
\begin{eqnarray*}
\|S(h)\|&\leq& \delta_0 t +K_1 R_{\delta_0 t} +K_2 R_{\delta_0 t} +\frac 32 R_{\delta_0 t}^2\\
&\leq& \delta_0 t +K_1 R_{\delta_0} +K_2 R_{\delta_0} +\frac 32 R_{\delta_0}^2\\
&<&\delta_0+K_1 R_{\delta_0} +K_2 R_{\delta_0} +\frac 32 R_{\delta_0}^2=2\,,
\end{eqnarray*}
so that $\|S(h)\|\leq 2s$ for some $s<1$. This shows (by the Ahlfors-Weill theorem) that the analytic part $h$ of every function $f=h+\overline g$ in the family $\mathcal F_{\delta_0 t}$ can be extended to a $K_s$-quasiconformal function in $\widehat\C$, where $K_s=(1+s)/(1-s)$. Using again that the family $\mathcal F_{\delta_0 t}$ is invariant under affine transformations, we get that not only $h$ but $h+ag$ (where $a\in\overline\D$) has a $K_s$-quasiconformal extension to $\widehat\C$. By arguing as in the proof of \cite[Thm. 2]{HM-QC}, we conclude that $f$ itself has a $K$-quasiconformal extension for an appropriate value of $K$.
\end{pf}


\begin{thebibliography}{10}

\bibitem{AW}
L. V. Ahlfors and G. Weill, A uniqueness theorem for Beltrami equations, \emph{Proc. Amer. Math. Soc.} \textbf{13} (1962), 975--978.

\bibitem{CHM}
M. Chuaqui, R. Hern\'andez, and M. J. Mart\'{\i}n, Affine and linear invariant families of harmonic mappings, arXiv:1405.5106 [math.CV].



\bibitem{Dur-Harm}
P. L. Duren, \emph{Harmonic Mappings in the Plane}, Cambridge University
Press, Cambridge, 2004.

\bibitem{GZ}
P. Ghatage and D. Zheng, Hyperbolic derivatives and generalized Schwarz-Pick estimates, \emph{Proc. Amer. Math. Soc.} \textbf{132} (2004), 3309--3318.



\bibitem{HM-Schw}
R. Hern\'andez and M. J. Mart\'{\i}n, Pre-Schwarzian and Schwarzian de\-ri\-va\-ti\-ves of harmonic mappings, \emph{J. Geom. Anal.} DOI 10.1007/s12220-013-9413-x. Published electronically on April 13th, 2013.

\bibitem{HM-QC}
R. Hern\'andez and M. J. Mart\'{\i}n, Quasi-conformal extersions of harmonic mappings in the plane, \emph{Ann. Acad. Sci. Fenn. Ser. A. I Math.} \textbf{38} (2013), 617--630.

\bibitem{HO}
T. Hosokawa and S. Ohno, Topological structures of the set of composition operators on
the Bloch space, \emph{J. Math. Anal. Appl.} \textbf{314} (2006), 736--748.


\bibitem{Lewy}
H. Lewy, On the non-vanishing of the Jacobian in certain one-to-one mappings, \emph{Bull. Amer. Math. Soc.} \textbf{42} (1936), 689--692.


\bibitem{Nehari}
Z. Nehari, The Schwarzian derivative and schlicht functions, \emph{Bull. Amer. Math. Soc.} \textbf{55} (1949), 545--551.

\bibitem{POM-I}
Ch. Pommerenke, Linear-invariante Familien analytischer Funktionen {I}, \emph{Math. Ann.} \textbf{155} (1964), 108--154.

\bibitem{POM-II}
Ch. Pommerenke, Linear-invariante Familien analytischer Funktionen {II}, \emph{Math. Ann.} \textbf{156} (1964), 226--262.


\bibitem{S-S}
T. Sheil-Small, Constants for planar harmonic mappings, \emph{J. London Math. Soc.} \textbf{42} (1990), 237--248.

\end{thebibliography}
\end{document}